\documentclass{article}
\usepackage{epsfig}
\usepackage{amsmath}
\usepackage{amsthm}
\usepackage{amssymb}
\usepackage{latexsym,amsfonts}
\usepackage{amscd}
\usepackage[all]{xy}
\parskip\smallskipamount
\numberwithin{equation}{section}

\newcommand{\kod}{\kappa}
\renewcommand{\mod}{\mathrm{mod}}

\newcommand{\mH}{H}

\newcommand{\mSpec}{\mbox{Spec\,}}

\newcommand{\Cc}{{\Bbb C}}
\newcommand{\Dd}{{\Bbb D}}
\newcommand{\Ee}{{\Bbb E}}

\newcommand{\Pp}{{\Bbb P}}
\newcommand{\Qq}{{\Bbb Q}}

\newcommand{\Xx}{{\Bbb X}}

\newcommand{\Zz}{{\Bbb Z}}

\newcommand{\ck}{{{\Oscr}}} 
 
\newcommand{\cA}{{\cal{A}}}
\newcommand{\cB}{{\cal{B}}}

\newcommand{\cD}{{\cal{D}}}

\newcommand{\cO}{{\cal{O}}}
\newcommand{\cP}{{\cal{P}}}

\newcommand{\lra}{\longrightarrow}

\newcommand{\ra}{\rightarrow}

\newtheorem{Prop}{Proposition}[section]
\newtheorem{Thm}[Prop]{Theorem}
\newtheorem{Lemma}[Prop]{Lemma}
\newtheorem{Cor}[Prop]{Corollary}
\newtheorem{Example}[Prop]{Example}
\newtheorem{Remark}[Prop]{Remark}
\newtheorem{Conjecture}[Prop]{Conjecture}

\let\cal\mathcal
\def\Ascr{{\cal A}}
\def\Bscr{{\cal B}}

\def\Dscr{{\cal D}}
\def\Escr{{\cal E}}
\def\Fscr{{\cal F}}

\def\Lscr{{\cal L}}

\def\Nscr{{\cal N}}
\def\Oscr{{\cal O}}
\def\Pscr{{\cal P}}

\def\Xscr{{\cal X}}

\let\blb\mathbb
\def\CC{{\blb C}}

\def\QQ{{\blb Q}}

\def \PP{{\blb P}}

\def \ZZ{{\blb Z}}

\def \NN{{\blb N}}
\def \RR{{\blb R}}

\def\Hom{\operatorname {Hom}}
\def\Ext{\operatorname {Ext}}
\def\RHom{\operatorname {RHom}}
\def\cone{\operatorname {cone}}
\def\coh{\operatorname {coh}}

\def\Aut{\operatorname {Aut}}
\def\Num{\operatorname {\Nscr}}

\def\Sl{\operatorname {Sl}}
\def\Pic{\operatorname {Pic}}
\def\Tor{\operatorname {Tor}}
\def\rad{\operatorname {rad}}
\def\rk{\operatorname {rk}}
\def\deg{\operatorname {deg}}

\def\ch{\operatorname {ch}}
\def\Td{\operatorname {Td}}
\def\Lotimes{\overset{L}{\otimes}}
\def\top{\operatorname{top}}
\def\Qcoh{\operatorname{Qcoh}}
\def\End{\operatorname{End}}
\def\Stab{\operatorname{Stab}}

\begin{document}
\title{Fourier-Mukai Transforms}
\author{Lutz Hille and Michel Van den Bergh}

\date{\today}
\maketitle
\begin{abstract}
 In this paper we discuss some of the recent developments on derived
equivalences in algebraic geometry. 
\end{abstract}
\tableofcontents
\section{Some background}\label{ref-1-0}
\label{ref-1-1}

In this paper we discuss some of the recent developments on derived equivalences
in algebraic geometry but
we don't intend to give any kind of comprehensive survey.
It is better to regard this paper as a set of pointers to
some of the recent literature. 

To put the subject in context we start with some historical
background.  Derived (and triangulated) categories were introduced by
Verdier in his thesis (see \cite{Deligne,Verdier}) in order to simplify
homological algebra. From this point of view the role of derived
categories is purely technical.

The first non-trivial derived equivalence  in the literature
is between the derived categories of sheaves on a sphere bundle and its dual bundle \cite{Sato}. The
equivalence resembles  Fourier-transform and is
now known as a ``Fourier-Sato'' transform.

The first purely algebro-geometric derived equivalence seems to appear
in \cite{Mukai} where is it is shown that an abelian variety $A$
and its dual $\hat{A}$ have equivalent derived categories of coherent sheaves. Again the equivalence is similar to a Fourier-transform and is
therefore called a ``Fourier-Mukai'' transform.

In 
\cite{Beilinson} Beilinson showed that 
$\PP^n$ is derived equivalent to  a (non-commutative) finite
dimensional algebra. This explained earlier results by Barth and Hulek
 on the relation between vector bundles and linear
algebra.  Beilinson's result has been generalized to other
varieties and has evolved into the theory of exceptional sequences (see
for example \cite{BondalOrlovsemiorthogonal}). 
The observation that derived equivalences do not preserve commutativity
is significant for non-commutative algebraic geometry (see for
example \cite{GeigleLenzing}).

Most algebraists probably became aware of the existence non-trivial
derived equivalences when Happel showed that ``tilting'' (as
introduced by Brenner and Butler \cite{BrennerButler}) leads to a
derived equivalence between finite dimensional algebras~\cite{Happel}. This was generalized by Rickard  who worked out the
Morita theory for derived categories of rings \cite{Rickard1,Rickard2}.

Hugely influential was the so-called homological mirror symmetry
conjecture by Kontsevich \cite{Kontsevichhomolmirror} which states
(very roughly) that for two Calabi-Yau manifolds $X$, $Y$ in a mirror
pair, the bounded derived category of coherent sheaves on $X$ is equivalent to
a certain triangulated category (the Fukaya category) related to the
symplectic geometry of $Y$. The homological mirror symmetry conjecture
was recently proved by Seidel for quartic surfaces (which are the
simplest Calabi-Yau manifolds after elliptic curves) \cite{Seidel}. 

Finally this introduction would be incomplete without at least mentioning the
celebrated Riemann-Hilbert correspondence 
\cite{borel,kashiwara,mebkhout2,mebkhout1}
which gives a derived
equivalence between  sheaves of vector spaces and regular holonomic D-modules on a complex
manifold or a smooth algebraic variety (depending on context). This is
a far reaching generalization of the classical correspondence between 
local systems and vector bundles with flat connections.

{\sl Acknowledgment. } The authors would like to thank Dan
Abramovich, Paul Balmer, Alexei Bondal, Tom Bridgeland, Daniel
Huybrechts, Pierre Schapira, Paul Smith and the anonymous referee for helpful comments on the
first version of this paper.

There are many other survey papers dedicated to Fourier-Mukai
transforms.  We refer in particular to Raphael Rouquier's
``Cat\'egories d\'eriv\'ees et g\'eometrie algebriques''
\cite{rouquier}. Another good source of information is given by preliminary
course notes by Daniel Huybrechts \cite{Huybrechts}.

\section{Notations and conventions}
\label{ref-2-2}
Throughout we work over the base field $\Cc$. The bounded
derived category of 
coherent sheaves on a variety $X$ is denoted by $\cD^b(X)$. Similarly, the
bounded derived 
category of finitely generated modules over an algebra $A$ is denoted by
$\cD^b(A)$. The shift functor in the derived 
category is denoted by $[1]$.  All functors between triangulated
categories are additive and exact (i.e.\ they commute with shift
and preserve distinguished triangles).

A sheaf is a coherent $\cO_X$--module
and a point in $X$ is always a closed point. The structure sheaf of a
point $x$ will be denoted by $\ck_x$. The canonical divisor of a
smooth projective variety is denoted by $K_X$ and the canonical sheaf is
denoted by $\omega_X$. 

\section{ Basics on Fourier-Mukai transforms}
\label{ref-3-3}

Let $X$ and $Y$ be connected smooth projective varieties. We are
interested in equivalences of the derived categories $\Phi: \cD^b(Y)
\lra \cD^b(X)$. Such varieties $X$ and $Y$ are also called {\sl
  Fourier-Mukai partners} and the equivalence $\Phi$ is called a {\sl
  Fourier-Mukai transform}. In this section we will discuss
some properties which remain invariant under Fourier-Mukai transforms.
The main technical tool is Orlov's theorem (see below) which states
that any derived equivalence  $\Phi: \cD^b(Y)
\lra \cD^b(X)$ is coming from a complex on the product $Y\times X$.

Given Fourier-Mukai $X,Y$ it is also interesting to precisely classify
the Fourier-Mukai transforms $\cD^b(Y) \lra \cD^b(X)$ (it is usually
sufficient to consider $X=Y$). This is generally
a much harder problem which has been solved in only a few special
cases, notably abelian varieties \cite{Orlovabelian} and varieties
with ample canonical or anti-canonical divisor (see Theorem
\ref{ref-4.4-14} below).

\medskip

To start one has the following simple result.
\begin{Lemma}[{\cite[Lemma 2.1]{BridgeMaciociasurfaces}}]
\label{ref-3.1-4}
If $X$ and $Y$ are Fourier-Mukai partners, then $\dim(X) = \dim(Y)$
and the canonical line bundles $\omega_X$ and $\omega_Y$ have the same
order.
\end{Lemma}
\begin{proof}
  The proof is an exercise in the use of \emph{Serre functors}
  \cite{BondalKapranov}. The Serre functor $S_X=-
  \otimes \omega_X[\dim(X)]$ on $X$  is uniquely characterized by the
  existence of natural isomorphisms
\begin{equation}
\label{ref-3.1-5}
\Hom_{\Dscr^b(X)}(\Escr,\Fscr)\cong \Hom_{\Dscr^b(X)}(\Fscr,S_X\Escr)^\ast
.
\end{equation}
By uniqueness it is clear that any  Fourier-Mukai transform 
commutes with Serre functors. Pick a point $y\in
Y$ and put $\Escr=\Phi(\Oscr_y)$. The fact that $S_Y[-\dim
Y](\Oscr_y)\cong \Oscr_y$ yields $S_X[-\dim Y](\Escr)\cong \Escr$, or
$\Escr\otimes_X \omega_X[\dim X-\dim Y] \cong \Escr$. Looking at the
homology of $\Escr$ we see that this impossible if $\dim Y\neq \dim
X$. The statement about the orders of $\omega_X$ and $\omega_Y$
follows by considering the orders of the functors $S_X[-\dim X]$ and
$S_Y[-\dim Y]$.
\end{proof}
The following important result tells that any derived equivalence
between $\cD^b(Y)$ and $\cD^b(X)$ is obtained from an object on the
product $Y\times X$.
\begin{Thm} [\cite{OrlovFMtransf}] \label{ref-3.2-6}
Let $\Phi: \cD^b(Y) \lra \cD^b(X)$ be a fully faithful functor. Then
there exists an object $\cP$ in $\cD^b(Y \times X)$, unique up to
isomorphism, such that $\Phi$ is isomorphic to the functor
$$
\Phi^{\cP}_{Y \ra X}(-) := \pi_{X \ast} (\cP \otimes_{\Oscr_{Y\times X}}
\pi_Y^{\ast}(-)), 
$$
where $\pi_X$ and $\pi_Y$ are the projection maps and $\pi_{X \ast}$,
$\otimes$, and $\pi_Y^{\ast}$ are the appropriate derived functors. 
\end{Thm}
In the original statement of this theorem $\Phi$ was required to have
a right adjoint but this condition is automatically fulfilled by
\cite{BondalVandenBergh,BondalKapranov}.

\medskip

The object $\cP$ in the theorem above is also called the \emph{kernel}
of the Fourier-Mukai transform.

\begin{Remark}
  Theorem \ref{ref-3.2-6} is quite remarkable as for example its analogue
  for affine varieties or finite dimensional algebras is unknown
  (except
for hereditary algebras \cite{MiyachiYekutieli}).
  Projectivity is used
  in the proof in the following way: let $\Lscr$ be an ample line
  bundle on a projective variety $X$. Then for any coherent sheaf
  $\Fscr$ on $X$ one has $\Hom_{\coh(X)}( \Fscr,\Lscr^{-n})=0$ for
  large~$n$. If $X$ is for example affine then $\Oscr_X$ is ample
but this additional property does not hold.
\end{Remark}

It would seem useful to generalize
Theorem \ref{ref-3.2-6} to singular varieties, in particular those
occurring in the minimal model program (see below). A first result in
this direction has been obtained by Kawamata \cite{Kawamata} who
proves the analogue of Theorem \ref{ref-3.2-6} for orbifolds.

\medskip

The real significance of Theorem \ref{ref-3.2-6} is that it makes it
possible to define $\Phi$ on objects functorially derived from $X$ and
$Y$.  For example (see \cite{Caldararu,OrlovRussianMathSurveys}) let
$\ch_X'(-)=\ch_X(-).\Td(X)^{1/2}$ (where $\ch_X(-)$ is the Chern
character and $\Td(X)$ is the Todd class of $X$). Using $\ch_{Y\times
  X}'(\Pscr)$ as kernel one finds a linear isomorphism of vector
spaces
$$
\mH^{\ast}(\Phi): \mH^{\ast}(Y,\Qq) \lra \mH^{\ast}(X,\Qq)
$$
preserving parity of degree. Since the Chern character of $\Pscr$
and the Todd class on $Y\times X$ may have denominators the same
result is not a priory true for $H^\ast(X,\ZZ)$.  However it is true
for elliptic curves (trivial)
and for abelian and K3-surfaces \cite{MukaiK3}.

\begin{Remark}
  In order to circumvent the non-preservation of integrality it may
  be convenient to replace $H^\ast(X,\ZZ)$ by topological K-theory
  \cite{Karoubi} $K^\ast(X)^{\top}=K^0(X)^{\top}\oplus
  K^1(X)^{\top}$ which is the K-theory of complex vector bundles (not
  necessarily holomorphic) on the underlying real manifold of $X$.
  Topological K-theory is a cohomology theory satisfying the usual
  Eilenberg-Steenrod axioms except the dimension axiom (which fixes
  the cohomology of a point).  Since $K^\ast(-)^{\top}$ has the appropriate
  functoriality properties \cite{Karoubi} one proves that
 $\Phi$ induces an isomorphism
\[
K^\ast(\Phi)^{\top}:K^\ast(Y)^{\top}\rightarrow K^\ast(X)^{\top}
\]
It
follows from the Atiyah-Hirzebruch spectral sequence 
  that $K^\ast(X)^{\top}$ is a finitely
generated $\ZZ/2\ZZ$ graded abelian group such that the
Chern-character
\[
\ch:K^\ast(X)^{\top}\rightarrow H^\ast(X,\QQ)
\] 
induces an isomorphism \cite[Eq (3.21)]{Hilton}
\[
K^\ast(X)^{\top}\otimes_\ZZ\QQ\cong H^\ast(X,\QQ)
\]
In good cases the lattices given by $K^\ast(X)^{\top}$ and
$H^\ast(X,\ZZ)$ are the same. This is for example the case
for curves, K3 surfaces and abelian varieties.
\end{Remark}

By Riemann-Roch the following diagram is commutative
\[
\begin{CD}
 K^0(Y) @>K^0(\Phi)>> K^0(X)\\
@VV\ch'_Y(-) V @VV\ch'_X(-) V \\
 \mH^{\ast}(Y,\Qq) @>\mH^{\ast}(\Phi)>> \mH^{\ast}(X,\Qq)
\end{CD}
\]
$K^0(X)$ is equipped with the so-called Euler form
\[
e([E],[F])=\sum_i (-)^i\dim \Hom_{\Dscr^b(X)} (E,F[i])
\]
which is of course preserved by $K^0(\Phi)$. 
The map $\ch'_X(-)$ is compatible with the Euler form up to sign
 provided one twist the standard bilinear form on cohomology
(obtained from Poincare duality) slightly \cite{Caldararu}. More precisely put
\[
\check{v}=i^{\deg v} e^{-(1/2) K_X} v
\]
and
\[
\langle v, w\rangle=\deg(\check{v}\cup w)
\]
Then 
\[
e([E],[F])=-\langle \ch'_X(E),\ch'_X(F)\rangle
\]
The map $H^i(\Phi)$ is an isometry for $\langle-,-\rangle$.

The standard grading on $H^\ast(X,\CC)$ is of course not preserved
by a Fourier-Mukai transform. However there is a different grading
with is preserved. Define
\[
{}^nH^\ast(X,\CC)=\bigoplus_{j-i=n} H^{i,j}(X)
\]
where $H^m(X,\CC)=\oplus_{i+j=m}
H^{i,j}(X,\CC)=\oplus_{i+j=m}H^i(X,\Omega_X^j)$ is the Hodge
decomposition \cite[\S0.6]{GriffithsHarris}. It is classical that algebraic
cycles lie in ${}^0H^\ast(X,\CC)$. From the fact that the kernel of
$\mH^{\ast}(\Phi)$ is algebraic it follows that $H^\ast(\Phi)$
preserves the ${}^\ast(-)$ grading.

\medskip

As another application of functoriality note that if $S$ is of finite type
then there is an equivalence 
\[
\Phi_S:\Dscr^b(Y_S)\rightarrow \Dscr^b(X_S)
\]
induced by $\Pscr_S$ (i.e.\ a Fourier-Mukai transform
extends to families).

\begin{Example} Here we give an example of a Fourier-Mukai transform
  which is very important for mirror-symmetry. Assume first that $Z$ is a
  four dimensional symplectic manifold and let $i:S^2\rightarrow Z$ be an
  embedding of a sphere as a Lagrangian submanifold. Then there exists
  a symplectic automorphism $\tau$ of $L$ which is trivial 
 outside a tubular neighborhood of $S^2$ 
  and which is the antipodal map on  $S^2$ itself
  \cite{Seidel1}. $\tau$ is called the \emph{symplectic Dehn twist}
  of $Z$ associated to~$i$.
  
  By the homological mirror symmetry conjecture there should be an
  analogous notion for derived categories of varieties. This was
  worked out in \cite{SeidelThomas} (see also \cite{KS,RZ}). It turns out that the analogue
  of a Lagrangian sphere is a so-called spherical object. To be more
  precise $\Escr\in \cD^b(X)$ is \emph{spherical} if
  $\Hom^i_{\cD^b(X)}(\Escr,\Escr)$ is equal to $\CC$ for $i=0,\dim X$
  and is zero in all other degrees and if in addition $\Escr\cong 
\Escr\otimes \omega_X$. 

Associated to a spherical object $\Escr\in \cD^b(X)$ there is an
auto-equivalence $T_\Escr$ of $\cD^b(X)$, informally defined by
\[
T_\Escr(\Fscr)=\cone\left(\RHom_{\cD^b(X)}(\Escr,\Fscr)\otimes_\CC \Escr
\xrightarrow{\text{\upshape evaluation}} \Fscr\right)
\]
The non-functoriality of cones leads to a slight technical problem with the
naturality of this
definition. This would
be a problem for abstract triangulated categories but it can be
rectified here using the fact that $\cD^b(X)$ (being a
derived category) is the $H^0$-category of an exact DG-category.

It is easy to show that the kernel of $T_\Escr$ is given by
\[
\cone\left(\check{\Escr}\boxtimes \Escr\xrightarrow{\phi} \Oscr_\Delta\right)
\]
where $\check{\Escr}=\RHom_{\Oscr_X}(\Escr,\Oscr_X)$, $\Oscr_\Delta$
is the structure sheaf of the diagonal and $\phi$ is the obvious map.

If $X$ is a K3-surface then $\Oscr_X$ is spherical and the kernel of
$T_{\Oscr_X}$ is given by $\Oscr_X(-\Delta)$. Other examples of
spherical objects are structure sheaves of a rational curve on a smooth surface with self
intersection $-2$ and restrictions of exceptional objects to
anticanonical divisors. In particular this last construction yields
spherical objects on hypersurfaces of degree $n+1$ in $\PP^n$.
\end{Example}

It is convenient to have a criterion for a functor of the form 
$\Phi^{\cP}_{Y \ra X}(-) := \pi_{X, \ast} (\cP \otimes
\pi_Y^{\ast}(-))$ to be an equivalence. The following result
originally due to Bondal and Orlov
\cite{BondalOrlovsemiorthogonal} and slightly amplified by Bridgeland
\cite[Theorem 1.1]{Bridgelandequivtriang} shows that we can use the skyscraper
sheaves as test objects.

\begin{Thm}\label{ref-3.6-7}
  Let $\cP$ be an object in $\cD^b(Y\times X)$. Then the functor
   $\Phi := \Phi^{\cP}_{Y \ra X}(-):
  \cD^b(Y) \lra \cD^b(X)$ is fully faithful if and only if  the following
conditions hold
\begin{enumerate}
  \item for each point $y$ in $Y$
$$
\Hom_{\cD^b(X)}(\Phi(\ck_y),\Phi(\ck_y)) = \Cc 
$$
\item for each pair of points $y_1$ and $y_2$ and each integer $i$
$$
\Hom_{\cD^b(X)}^i(\Phi(\ck_{y_1}),\Phi(\ck_{y_2})) = 0 \mbox{ unless
  } y_1 = y_2 \mbox{ and } 0 \leq i \leq \dim Y.
$$
\end{enumerate}
If these conditions hold then $\Phi$ is an equivalence if and only if 
 $\Phi(\ck_y) \otimes \omega_X \cong \Phi(\ck_y)$ for all $y\in Y$.
\end{Thm}
\begin{Remark} Assume that $\Pscr$ is an object in $\coh(Y\times X)$
flat over $Y$ and write
$\Pscr_y=\Phi(\Oscr_y)$. Then the previous theorem implies
that $\Phi$ is fully faithful if and only if
\begin{enumerate}
  \item for each point $y$ in $Y$
$$
\Hom_{\cD^b(X)}(\Pscr_y,\Pscr_y) = \Cc 
$$
\item for each pair of points $y_1\neq y_2$ and each integer $i$
$$
\Ext^i_{\Oscr_X}(\Pscr_{y_1},\Pscr_{y_2}) = 0.
$$
\end{enumerate}
\end{Remark}

It is obvious that the conditions for Theorem \ref{ref-3.6-7} are necessary.
Proving that they are also sufficient is much harder. Since the proof
in \cite{BondalOrlovsemiorthogonal} only works for
derived categories of coherent sheaves, we make explicit some
of the steps in Bridgeland's proof (see \cite{Bridgelandequivtriang})
which are valid for more general triangulated categories.

Let $\cA$ be a triangulated category. A subset $\Omega$ is called \emph{spanning} if for each object $a$ in $\cA$ each of the following
conditions implies $a=0$: 
\begin{enumerate}
\item $\Hom^i(a,b) = 0$ for all $b \in \Omega$ and all $i \in \Zz$,
\item $\Hom^i(b,a) = 0$ for all $b \in \Omega$ and all $i \in \Zz$.
\end{enumerate}
It is easy to see that the set of all skyscraper sheaves on a smooth
projective variety $X$ is a spanning class for $\cD^b(X)$.  Note
that a spanning class will not usually generate $\Ascr$
in any reasonable sense. 
\begin{Thm}[{\cite[Theorem 2.3]{Bridgelandequivtriang}}]
\label{ref-3.8-8}
Let $F: \cA \lra \cB$ be an exact functor between triangulated
categories with left and right adjoint. Then $F$ is fully faithful if
and only if there exists a spanning class $\Omega$ for $\cA$ such that
for all elements $a_1, a_2$ in $\Omega$, and all integers $i$, the
homomorphism
$$
F: \Hom_{\cA}^i(a_1, a_2) \lra \Hom_{\cB}^i(Fa_1, Fa_2)
$$
is an isomorphism.
\end{Thm}

Recall that a category is called indecomposable if it is not
the direct sum of two non-trivial subcategories. The derived category
$\cD^b(X)$ is indecomposable for $X$ connected. For a finite
dimensional algebra $A$ the derived category $\cD^b(A)$ is connected
precisely when $A$ is connected.

\begin{Thm}[{\cite[Theorem 2.3]{BKR}}]
\label{ref-3.9-9} Let $F:\Ascr\rightarrow \Bscr$ be a fully
  faithful functor between triangulated categories with Serre functors
  $S_\Ascr$, $S_\Bscr$ (see \eqref{ref-3.1-5}) possessing a left adjoint.
  Suppose that $\Ascr$ is non-trivial and $\Bscr$ is indecomposable.
  Let $\Omega$ be a spanning class for $\Ascr$ and assume that 
$FS_\Ascr(\omega)
\cong S_\Bscr F(\omega)$ for all $\omega\in \Omega$. Then
$F$ is a equivalence of categories.
\end{Thm}

It follows from \cite{BondalVandenBergh,BondalKapranov} that $\Phi^\Pscr_{Y\rightarrow
  X}$ has both a right and a left adjoint. Explicit formulas are
for the left and the right adjoint are \cite[Lemma 4.5]{Bridgelandequivtriang}:
$$
\Phi^{\check\cP \otimes \pi^{\ast}_X \omega_X[\dim X]}_{X \ra Y}(-)
\mbox{\ \ and\ \ }
 \Phi^{\check \cP \otimes \pi^{\ast}_Y \omega_Y[\dim Y]}_{X \ra Y}(-)
$$

Applying Theorems \ref{ref-3.8-8},\ref{ref-3.9-9} with
$F=\Phi^{\Pscr}_{Y\rightarrow X}$ and $\Omega=\{\Oscr_y\mid y\in Y\}$
almost proves Theorem \ref{ref-3.6-7} except that we seem to need additional
information on $\Hom_{\cD^b(X)}^i(\Phi(\Oscr_y),\Phi(\Oscr_y))$ for
$i>0$.  It is not at all obvious but it turns out that this extra
information is unnecessary. Although it is not clear how to formalize
it, it seems that this part of the proof may generalize whenever $Y$
is the solution of some type of moduli problem in a triangulated
category $\Bscr$ (with $\Pscr$ being the universal family). See
\cite{Bridgelandflops,BKR,VandenBerghCrepant} for other manifestations of this
principle.

\section{The reconstruction theorem}
\label{ref-4-10}
It is quite trivial to reconstruct $X$ from the abelian category
$\coh(X)$ \cite{gabriel,rosenberg,rouquier}. For example the points of $X$ are in one-one correspondence
with the objects  in $\coh(X)$ without proper subobjects. With a little
more work one can also recover the Zariski topology on $X$ as well as
the structure sheaf.

It is similarly  of
interest to know to which extent one can reconstruct a variety from its
derived category. The existence of non-isomorphic Fourier-Mukai partners
shows that this cannot be done in general, but it is possible
 if  the canonical sheaf or the anticanonical sheaf is ample. Later
 Balmer and Rouquier \cite{Balmer,rouquier} have shown independently that one can reconstruct the variety from the
 category of coherent sheaves viewed as a tensor category, the crucial
 point is that the tensor product allows to reconstruct the point
 objects for any variety.

\begin{Thm}[{\cite[Theorem 2.5]{BondalOrlovreconstruction}}]
  \label{ref-4.1-11} 
  Let $X$ be a smooth connected projective variety with either
  $\omega_X$ ample or $\omega_X^{-1}$ ample. Assume $\cD^b(X)$ is
  equivalent to $\cD^b(Y)$. Then $X$ is isomorphic to $Y$.
\end{Thm}

\begin{proof}  We give a proof based on Orlov's theorem. For another
proof see \cite{rouquier}. Note that
$Y$ is also connected since $\cD^b(Y)\cong \cD^b(X)$ is connected. 

Let $\Phi:\Dscr^b(Y)\rightarrow \Dscr^b(X)$ be the derived equivalence
and let $S$ be the Serre functor $- \otimes
\omega_X[\dim X]$ on $X$. Recall that it is intrinsically defined
by \eqref{ref-3.1-5}.
We say that $E$ in $\cD^b(X)$ is a \emph{point object} if 
\begin{enumerate}
\item $E \cong S(E)[i]$ for some integer $i$, 
\item $\Hom^i(E,E) = 0$ for all $i <0$, and 
\item $\Hom(E,E) = \Cc$. 
\end{enumerate}
It is easy to prove that the only point objects in $\cD^b(X)$ (under
the assumptions on $\omega_X)$ are the shifts of the skyscraper
sheaves. The main point is 1., since this condition and the ampleness of $\omega_X^{\pm
  1}$ easily implies that $E$ has finite length cohomology.

It follows that $\Phi$ sends skyscraper sheaves to shifts
of skyscraper sheaves. Then the proof may then be finished using Corollary
\ref{ref-4.3-13}
below.
\end{proof}

We need the following standard fact.
\begin{Prop} \label{ref-4.2-12}
  Let $\pi: Z\rightarrow S$ be a flat morphism of schemes of finite type with $S$
  connected. Let $\Pscr\in \Dscr^-(\coh(Z))$ and and assume that for
  all $s\in S$ we have that $\Pscr\Lotimes_{\Oscr_Z}
  \pi^\ast\Oscr_s\cong \Oscr_z[n]$ for some $n\in \ZZ$, $z\in Z$. Then
  $\Pscr\cong i_\ast\Lscr[m]$ where $i:S\rightarrow Z$ is a section of $\pi$,
  $\Lscr\in \Pic(S)$ and $m\in\ZZ$.
\end{Prop}
\begin{proof}
  We claim first that the support of the cohomology $\Pscr$ is finite
  over $S$. Assume that this is false and let $H^i(\Pscr)$ be the
  highest cohomology group with non-finite support. Then, up to finite
  length sheaves we have $H^i(\Pscr)\otimes_{\Oscr_Z}
  \pi^\ast\Oscr_s\cong H^i(\Pscr\Lotimes_{\Oscr_Z} \pi^\ast\Oscr_s)$.
  Hence $H^i(\Pscr)\otimes_{\Oscr_Z} \pi^\ast\Oscr_s$ has finite
  length for all $s$ which is a contradiction.

It is now sufficient to prove that $\Pscr_0=\pi_\ast(\Pscr)$ is a shifted
line bundle given that $\Pscr_0\Lotimes_{\Oscr_S} \Oscr_s$ has
one-dimensional cohomology for all $s$. 

Fix $s\in S$ and assume
$\Pscr_0\Lotimes_{\Oscr_S} \Oscr_s\cong \Oscr_s[n]$. Using Nakayama's
lemma we deduce that there is a neighborhood $U$ of $s$ such
that $H^i(\Pscr_0\mid U)=0$ for $i>-n$.  We temporarily replace $S$ by
$U$.

 Applying $-\Lotimes_{\Oscr_S}
\Oscr_s$ to the triangle
\[
\tau_{\le -n-1} \Pscr_0\rightarrow \Pscr_0\rightarrow H^{-n}(\Pscr_0)[n]\rightarrow 
\]
we find $H^{-n}(\Pscr_0)\otimes_{\Oscr_S} \Oscr_s\cong \Oscr_s$ and
$\Tor_1^{\Oscr_S}(H^{-n}(\Pscr_0),\Oscr_s)=0$. Hence $H^{-n}(\Pscr)$ is a line
bundle on a neighborhood of $s$. Shrinking $S$ further we may assume
$\Pscr_0\cong \tau_{\le -n-1} \Pscr_0\oplus H^{-n}(\Pscr_0)[n]$ and hence
$\tau_{\le -n-1} \Pscr_0\Lotimes_{\Oscr_S} \Oscr_s=0$.   Shrinking
$S$ once again we have $\tau_{\le -n-1} \Pscr_0=0$ and thus $\Pscr_0\cong
H^0(\Pscr_0)[n]$ is a line bundle on a neighborhood of $s$.

Since this works for any $s$ and $S$ is connected we easily deduce
that $\Pscr_0$ is itself a shifted line bundle.
\end{proof}
We deduce the following 
\begin{Cor} \label{ref-4.3-13} Assume that $\Phi:\Dscr^b(Y)\rightarrow \Dscr^b(X)$ is a 
Fourier-Mukai transform between smooth connected projective varieties
which sends skyscraper sheaves to shifted skyscraper sheaves. Then 
$\Phi$ is of the form $\sigma_\ast(-\otimes_{\Oscr_X}\Lscr)[n]$ for
an isomorphism $\sigma:Y\rightarrow X$, $\Lscr\in \Pic(Y)$ and $n\in\ZZ$.
\end{Cor}
\begin{proof} By  Proposition \ref{ref-4.2-12} the kernel of $\Phi$
must be of the form $\Pscr=(1,\sigma_\ast)_\ast\Lscr[n]$ for some
map $\sigma:Y\rightarrow X$. The resulting $\Phi^\Pscr_{Y\rightarrow X}=
\sigma_\ast(-\otimes_{\Oscr_X})[n]$
will be a derived equivalence if and only if $\sigma$ is an isomorphism.
\end{proof}
One also obtains as a corollary the following result.
\begin{Thm}[{\cite[Theorem 3.1]{BondalOrlovreconstruction}}]
  \label{ref-4.4-14} 
Let $X$ be a smooth connected projective variety with ample canonical
or anticanonical sheaf. Then the group of isomorphism classes of
auto-equivalences of $\cD^b(X)$ is generated by the automorphisms of $X$,
the twists by line bundles and the  translations.
\end{Thm}

\begin{Remark}
It is clear that the notion of point object 
make sense for arbitrary triangulated categories with Serre functor. 

Let $\cD$ be the bounded derived category of modules over a
connected finite dimensional hereditary $\Cc$--algebra $A$. Then point
objects only exist 
for $A$ tame (or in the trivial case $A \cong \Cc$). For Dynkin
quivers or wild quivers the structure of the Auslander-Reiten
components is well-known (Gabriels work on Dynkin-quivers and Ringels
work on wild hereditary alegebras), consequently, point objects do not exist.
In the tame case the point objects 
are the shifts of quasi-simple modules in homogeneous tubes (see \cite{Ringeltame}). 
Let $A$ be not necessary hereditary and we assume $\cD^b(A)$ is
equivalent to $\cD^b(X)$ for some smooth projective variety $X$. Then
$\cD^b(A)$ has point objects.
The situation is
similar if we replace $X$ by a weighted projective variety. However,
it is an open 
problem to construct algebras $A$ having (sufficiently many) point
objects 
without knowing such an equivalence between $\cD^b(A)$ and
$\cD^b(X)$ for some (weighted) projective variety~$X$.
\end{Remark}
Note that there is a subtle point in the statement of Theorem
\ref{ref-4.1-11}. One does not \emph{apriori} require  $Y$ to have ample
canonical or anti-canonical divisor. If we preimpose this condition
then Theorem
\ref{ref-4.1-11} also follows from Theorem \ref{ref-4.6-15} below which
morally corresponds to the fact that derived equivalences
commute with Serre functors.

\begin{Thm}[\cite{OrlovRussianMathSurveys}]\label{ref-4.6-15} 
Let $X$ be a smooth projective variety. Then the integers
  $\dim\Gamma(X,\omega_X^{\otimes m})$ as well as the the
  canonical and anti-canonical rings are derived invariants.
\end{Thm}
Assume that $X$ is connected. For a Cartier divisor $D$ denote by
$R(X,D)$ the ring
\[
R(X,D)=\bigoplus_{n\ge 0} \Gamma(X,\Oscr_X(nD))
\]
and by $K(X,D)$ the part of degree zero of the graded quotient field
of $R(X,D)$. We have $K(X,D)\subset K(X)$ where $K(X)$ is the function
field of $X$. By \cite[Prop 5.7]{Ueno} $K(X,D)$ is algebraically
closed in $K(X)$. If some positive multiple of $D$ is effective then the
$D$-Kodaira dimension $\kappa(X,D)$ of $K(X,D)$ is the transcendence degree of
$K(X,D)$, otherwise we set $k(X,D)=-\infty$. It is clear that we have
\[
\kappa(X,D)\le \dim X
\]
and in case of equality we have $K(X)=K(X,D)$. 

The Kodaira dimension $\kappa(X)$ of $X$ is $\kappa(X,K_X)$.  $X$ is of general
type if
$\kappa(X,K_X)=\dim X$. 
\begin{Cor}[{\cite[Theorem 2.3]{KawamataKequivalence}}] \label{ref-4.7-16}
The Kodaira dimension is 
  invariant under Fourier-Mukai transforms. If $X$ is of general type
  then any Fourier-Mukai partner of $X$ is birational to $X$.
\end{Cor}
\begin{proof} This follows directly from  Theorem \ref{ref-4.6-15} and the
  preceding
discussion.
\end{proof}

\section{Curves and surfaces}
\label{ref-5-17}

In this section we consider Fourier-Mukai transforms for smooth
projective curves and smooth projective surfaces. For curves
the situation is rather trivial: only elliptic curves admit
non-trivial Fourier-Mukai transforms $\cD^b(C) \cong \cD^b(D)$, and
in that case the curves $C$ and $D$ must be isomorphic.  The group of
auto-equivalences of $\cD^b(C)$ is generated by the trivial ones and
the classical Fourier-Mukai transform (which is almost the same as the
auto-equivalence associated to the spherical object $\Oscr_E$).

For surfaces the
situation is more complicated and is worked out in detail
in \cite{BridgeMaciociasurfaces}. The classification of possible
non-trivial Fourier-Mukai transforms is based on the classification of
complex surfaces (see \cite[page 188]{BarthPetersVandeVen}).
This classification is summarized in Table 1. 
\medskip

\begin{figure}

\medskip

$$
\begin{array}{|l|c|c|c|c|c|}
\hline
\mbox{ Class of } \vrule width 0pt height 2.5ex depth 0.5ex X & \kod(X) & n_X & b_1(X)  & c_1^2 & c_2  \\
\hline \hline
\mbox{1) \vrule width 0pt height 2.5ex minimal rational }  & & & & & \\
\mbox{surfaces } &
                      -\infty &     & 0        & 8,9   & 4,3   \\        
\mbox{3) ruled surfaces }   & & & & & \\
\mbox{ of genus } g \geq 1 & 
                    -\infty &     & 2g         & 8(1-g)   & 4(1-g)   \\ 

\hline
\mbox{4)\vrule width 0pt height 2.5ex  Enriques surfaces } &
                    0 & 2    & 0        & 0   & 12   \\ 

\mbox{5) hyperelliptic surfaces } &
                   0 & 2,3,4,6         & 2   & 0   & 0   \\ 
\mbox{7) K3-surfaces } &
                   0 & 1 & 0  & 0 & 24 \\
\mbox{8) tori } &
                   0 & 1 & 4 & 0 & 0 \\
\hline
\mbox{9) \vrule width 0pt height 2.5ex minimal properly }  & & & & & \\
\mbox{elliptic surfaces} &
                   1 &   &   & 0 &  \geq 0 \\
\hline
\mbox{10) \vrule width 0pt height 2.5ex minimal surfaces }  & & & & & \\
\mbox{of general type } &
                   2    &  &\equiv 0\ \mod\ 2  & > 0 & >0 \\
\hline
\end{array}
$$
\centerline{ {\bf Table 1.} Classification of algebraic smooth complex surfaces}

\end{figure}

Let us start with the case of curves. 
Let $C$ be a smooth projective curve and denote by $g_C$ the genus of
$C$. According to the degree of the 
canonical divisor $K_C$ there are three distinct classes: 
\begin{enumerate}
\item $K_C < 0$: $C$ is the projective line $\Pp^1(\Cc)$ and $g_C = 0$, 
\item $K_C = 0$: $C$ is an elliptic curve and $g_C =1$, 
\item $K_C > 0$: $C$ is a curve of general type and $g_C > 1$.
\end{enumerate}

Using the reconstruction theorems  \ref{ref-4.1-11} and \ref{ref-4.4-14}
it is obvious that
non-trivial Fourier-Mukai transforms can only exist for elliptic
curves since $K_C^{-1}$ is ample in case 1.\ and $K_C$ is ample in case
3..

We will now look in somewhat more detail at the interesting case of
elliptic curves. Note that if $C$, $D$ are abelian varieties then it
is known precisely when $C$ and $D$ are derived equivalent and
furthermore the group $\Aut(\cD^b(C))$ consisting of auto-equivalences
of $\cD^b(C)$ (up to isomorphism) is also completely understood
\cite{Orlovabelian}. Here we give an elementary
account of the one-dimensional case. This is well-known and was
explained to us by Tom Bridgeland.
First we have the following result.
\begin{Thm} 
\label{ref-5.1-18} If $C$, $D$ are derived equivalent elliptic curves
then $C\cong D$.
\end{Thm}
\begin{proof} By the discussion in \S \ref{ref-3-3} the Hodge structures
on $H^1(C,\CC)$ and $H^1(D,\CC)$ are isomorphic. Since the
isomorphism
class of an elliptic curve is encoded in its Hodge structure on
$H^1(-,\CC)$
we are done.
\end{proof}
Determining the structure of $\Aut(\Dscr^b(C))$ requires slightly more
work.  For an elliptic curve $C$ let $e_C$ be the Euler form on
$K^0(C)$. By Serre duality $e_C$ is skew symmetric. Put
$\Num(C)=K^0(C)/\rad e_C\cong \ZZ^2$. $e_C$ defines a non-degenerate
skew symmetric form (i.e.\ a symplectic form) on $\Num(C)$ which we
denote by the same symbol.

$\Num(C)$ has
a canonical basis given by $v_1=[\Oscr_C]$, $v_2=[\Oscr_x]$ ($x\in C$ 
arbitrary).
The matrix of $e_C(v_i,v_j)_{ij}$ with respect to this basis is
\[
\begin{pmatrix}
0& 1\\
-1 & 0
\end{pmatrix}
\]
With respect to the 
standard basis the group of symplectic automorphisms of $\Num(C)$ may be
identified with $\Sl_2(\ZZ)$.

Let $T_1$, $T_2$ be the auto-equivalences of $C$ associated to the
spherical objects $\Oscr_C$ and $\Oscr_x$. It is not hard to see
that $T_2=-\otimes_{\Oscr_C} \Oscr_C(x)$ so only $T_1$ is a non-trivial
Fourier-Mukai transform.

One computes that with
respect to the standard basis the action of $T_1$, $T_2$ 
on $\Num(C)$ is given by matrices
\[
T_1=
\begin{pmatrix}
1&-1\\
0&1
\end{pmatrix}
\]
\[
T_2=
\begin{pmatrix}
1&0\\
1&1
\end{pmatrix}
\]
These matrices are standard generators for $\Sl_2(\ZZ)$
which satisfy the braid relation
\begin{equation}
\label{ref-5.1-19}
T_1T_2T_1=T_2T_1T_2
\end{equation}
\begin{Remark}
Since the objects $\Oscr_C$, $\Oscr_x$ form a
so-called $A_2$ configuration \cite{SeidelThomas} the
relation \eqref{ref-5.1-19} actually holds in $\Aut(\Dscr^b(C))$.
\end{Remark}
We have:
\begin{Thm} 
Let $\Aut^0(\Dscr^b(C))$ be the subgroup of $\Aut(\Dscr^b(C))$ consisting
of auto-equivalences of the form $\sigma_\ast(-\otimes_{\Oscr_C}\Lscr)[n]$
where $\sigma\in \Aut(C)$, $\Lscr\in\Pic^0(C)$ and $n\in 2\ZZ$. Then
the symplectic action of $\Aut(\Dscr^b(C))$ on $\Num(C)$ yields
an exact sequence
\[
0\rightarrow \Aut^0(\Dscr^b(C))\rightarrow \Aut(\Dscr^b(C))
\rightarrow \Sl_2(\ZZ)\rightarrow 0.
\]
\end{Thm}
\begin{proof} The existence of $T_1$, $T_2$ implies that the map  
  $\Aut(\Dscr^b(C)) \rightarrow \Sl_2(\ZZ)$ is onto.

Assume that $\Phi\in \Aut(\Dscr(C))$ act trivially on $\Num(C)$.
It is easy to see that
for an object $\Escr\in \Dscr^b(C)$ this implies
\begin{equation}
\label{ref-5.2-20}
\begin{split}
\deg \Phi(\Escr)&=\deg \Escr\\
\rk \Phi(\Escr) &=\rk \Escr
\end{split}
\end{equation}
The abelian category $\coh(D)$ is hereditary and hence every object in
$\Dscr^b(D)$ is the direct sum of its cohomology. Since
$\Phi(\Oscr_y)$ must be indecomposable we deduce from
\eqref{ref-5.2-20} that $\Phi(\Oscr_y)$ is a twisted skyscraper sheaf.

We find by Corollary \ref{ref-4.3-13} that $\Phi=
\sigma_\ast(-\otimes_{\Oscr_C}\Lscr)[n]$. The fact that $\Phi$ acts
trivially on $\Num(C)$ implies $\deg\Lscr=0$ and $n$ is even.
\end{proof}

\begin{Remark} Using similar arguments as above it is easy to see
  that the orbits of the action $\Aut(\Dscr^b(C))$ on the
  indecomposable objects in $\Dscr^b(C)$ are indexed by $\NN\setminus \{0\}$. The quotient
map is given by
\[
E\mapsto \gcd(\rk(E),\deg(E))
\]
In particular any indecomposable vector bundle is in the orbit
of an indecomposable finite length sheaf.
\end{Remark}

\begin{Remark}
The situation for elliptic curves is 
very similar to the situation for tubular algebras
\cite{Ringeltame,RingelHappelTubular}, tubular canonical algebras,  or
tubular weighted projective curves (weighted projective curves of
genus one)
\cite{MeltzerLenzingwtprojlinegenus1}. We quickly explain how these
three categories 
$\cD^b(C)$ ($C$ an elliptic curve), $\cD^b(\Xx)$ ($\Xx$ a tubular
weighted projective curve) and $\cD^b(\Lambda)$ ($\Lambda$ a tubular
canonical algebra or a tubular algebra) are related to each other. Any
elliptic curve $C$ 
admits a non-trivial automorphism $\phi:C \lra C$ $x \mapsto -x$. Let
$G\cong\Zz/2\ZZ$, generated by $\phi$. The category of
$G$-equivariant sheaves on $C$ is isomorphic to the category of
coherent sheaves on a weighted projective line of type $\widetilde
\Dd_4$. For the remaining types $\Ee_{6,7,8}$ we consider elliptic
curves with complex multiplication of order $3$, $4$ or $6$,
respectively. Then an analogous result holds for those curves (see
also \cite{Schiffmann}).
\end{Remark}

Now we discuss the case of surfaces.  In the rest of this section a
surface will be a smooth projective  surface.

Remember that a surface $X$ is called minimal if it
does not contain an exceptional curve $C$ (i.e.\ a smooth rational
curve with self intersection $-1$).  The possible non-trivial
Fourier-Mukai partners for minimal surfaces were classified by
Bridgeland and Maciocia in \cite{BridgeMaciociasurfaces}. This
classification is based on the classification of  surfaces (see
\cite[page 188]{BarthPetersVandeVen}) as summarized in Table 1
(we have only listed the algebraic surfaces as these are the only
ones of interest to us).

Table 1 is in
terms
of some standard invariants which we first describe. We have already
mentioned the Kodaira
dimension $\kod(X)$. It is either $-\infty,0,1$ or $2$ and divides the
minimal surfaces into four classes.  For an arbitrary surface $X$ there
is always a map $X\rightarrow X_0$ to a minimal surface. If $k(X)\ge 0$ then
$X_0$ depends only on the birational equivalence class of $X$
\cite[Proposition (4.6)]{BarthPetersVandeVen}. 

Further invariants are the first Betti
number $b_1(X)=\dim H^1(X,\CC)$, the square of the first Chern class 
$c_1^2(X)=K^2_X$ and the
second Chern class $c_2(X)$ (where $c_i=c_i(T_X)$). 
Finally, for surfaces
of Kodaira dimension zero one also needs the smallest natural number $n_X$
with $n_X K_X = 0$.

The invariants $b_1(X)$, $c_1(X)^2$, $c_2(X)$ contain exactly the
same information as the (numeric) Hodge diamond of $X$:
\[
\begin{array}{ccccc}
&&1&&\\
&q(X)&& q(X)&\\
p_g(X) && h^{1,1}(X) && p_g(X)\\
&q(X)&& q(X)&\\
&&1&&
\end{array}
\]

where $p_g(X)$ is the geometric genus of $X$, $q(X)$ is the
Noether number of $X$ and $h^{ij}(X)=\dim H^{ij}(X,\CC)$. One has
\begin{align*}
b_1(X)&=2q(x)\\
c_2(X)&=2+2p_g(X)-4q(X)+h^{1,1}(X)\\
\frac{1}{12}(c_1(X)^2+c_2(X))&=1-q(x)+p_g(X))
\end{align*}
The second line is  the Gauss-Bonnet
formula \cite[\S3.3]{GriffithsHarris} which says that $c_2(X)$ is equal
to the Euler number $\sum_i \dim (-1)^i \dim H^i(X,\CC)$ of $X$. The
third formula is Noether's formula. It follows from applying the
Riemann-Roch theorem \cite[Thm I.(5.3)]{BarthPetersVandeVen} to the structure sheaf.

For abelian and K3-surfaces the so-called transcendental lattice
is of interest.  First note that $H^2(X,\ZZ)$ is free. For abelian
surfaces this is clear since they are tori and for K3 surfaces it
is \cite[Prop VIII(3.2)]{BarthPetersVandeVen}.
The Neron-Severi lattice is 
$
N_X=H^2(X,\ZZ)\cap H^{1,1}(X)
$
and the transcendental lattice $T_X$ is the sublattice of $H^2(X,\ZZ)$
orthogonal to $N_X$.

\begin{Thm}[{\cite[Theorem 1.1]{BridgeMaciociasurfaces}}] ~ 
\label{ref-5.6-21}
  Let $X$ and $Y$ be a non-isomorphic smooth connected complex
  projective surfaces with equivalent derived categories $\cD^b(X)$
  and
  $\cD^b(Y)$ such that $X$ is minimal. Then either 
\begin{enumerate}
  \item $X$ is a torus (an abelian surface, in class 8)) and $Y$ is also
  a
  torus with Hodge-isometric transcendental lattice, 
  \item $X$ is a K3-surface (a surface in class 7)) and $Y$ is also
  a K3-surface with Hodge isometric transcendental lattice, or 
\item
  $X$ is an elliptic surface and $Y$ is another elliptic surface
  obtained by taking a relative Picard scheme of the elliptic
  fibration on $X$.
\end{enumerate}
\end{Thm}
A Hodge isometry between transcendental lattices is an isometry 
under which the
one dimensional subspaces $H^0(X,\omega_X)$ and $H^0(Y,\omega_Y)$
of $T_X\otimes_\RR \CC$ and $T_Y\otimes_\RR\CC$ correspond.

\medskip

The proof of Theorem \ref{ref-5.6-21} is quite involved and uses case
by case analysis quite essentially.  As a very rough indication of
some of the methods one might use, let us show that
if $X$ is minimal then so is $Y$ and $X$ and $Y$ are in the same class. Along
the way we will settle the easy case $\kod(X)=2$.

\medskip

{\sl Step 1: } By Corollary \ref{ref-4.7-16} and the discussion in \S \ref{ref-3-3} $X$
  and $Y$ have the same Kodaira dimension and the same Hodge diamond.
  In particular they have the same $b_1(-)$, $c_1(-)^2$ and
  $c_2(-)$. Hence
if they are both minimal then they are in the same class.

{\sl Step 2: }
 Assume now that $X$ is minimal and let $Y\rightarrow Y_0$ be a
  minimal
model of $Y$. We have $b_1(Y)=b_1(Y_0)$ \cite[Theorem I.(9.1)]{BarthPetersVandeVen}.
If $\kod(X)=-\infty,1,2$ then the class of $X$ is recognizable from
$b_1(X)$ and hence $Y_0$ must be in the same class as $X$. If $Y_0$ is not
in class 1,10) then  it follows from the classification that $c_1(Y_0)^2=
c(X)^2$ and hence $c_1(Y_0)^2=c_1(Y)^2$.  If $Y_0$ is  in class 10) then 
by Corollary \ref{ref-4.7-16} we have $X=Y_0$ and hence we also have $c_1(Y_0)^2=c_1(Y)^2$.
Since $c_1(-)^2$ changes by one under a blowup
\cite[Theorem I.(9.1)(vii)]{BarthPetersVandeVen} it follows in these cases that
$Y=Y_0$. 

If $Y_0$ is is in class 1) then in principle we could have $c_1(Y_0)^2=9$,
$c_1(Y)^2=c_1(X)^2=8$. But then in $Y$ is the blowup of $\PP^2$ in a
point and hence is Del-Pezzo. We conclude by the reconstruction
 theorem \ref{ref-4.1-11} that $X=Y$ which is a contradiction.

{\sl Step 3: }
 If $\kod(X)=0$ then $\omega_X$ has finite order and hence
the same is true for $Y$ by Lemma \ref{ref-3.1-4}. This is impossible if $Y$
is not minimal. 

\medskip

Let us also say a bit more on the K3 and abelian case. Assume that
$X$ is a a K3 or abelian surface. Then according \cite{MukaiK3} the
Chern character $K^0(X)\rightarrow H^\ast(X,\QQ)$ takes it values in
$H^\ast(X,\ZZ)$. As before let $\Num(X)$ be $K^0(X)$ modulo the
radical of the Euler form. Since the intersection form on
$H^\ast(X,\ZZ)$ is non-degenerate it follows that $\Num(X)$ is
the image of $K^0(X)$ in $H^\ast(X,\ZZ)$. It is easy to see that
the orthogonal to $\Num(X)$ is $T_X$.

Now assume that $X$ and $Y$ are derived equivalent K3 or abelian
surfaces.  Again by \cite{MukaiK3} the induced isometry between 
$H^\ast(X,\QQ)$ and $H^\ast(Y,\QQ)$ yields an isometry between
$H^\ast(X,\ZZ)$ and $H^\ast(X,\ZZ)$.
By the above discussion there is an isometry between
$T_X$ and $T_Y$. This is a Hodge isometry since
$H^0(X,\omega_X)={}^{2} H^\ast(X,\CC)$.

The complete result for K3 or abelian surfaces is as follows.
\begin{Thm}[\cite{OrlovFMtransf}, see also \cite{BridgeMaciociasurfaces}]
Let $X$ and $Y$ be a pair of either K3-surfaces or abelian surfaces
(tori) then the following statements are equivalent. 
\begin{enumerate}
 \item There exists a Fourier-Mukai transform $\Phi: \cD^b(Y)\lra
 \cD^b(X)$. 
 \item There is an Hodge isometry $\phi^t: T(Y) \lra T(X)$.
 \item There is an Hodge isometry $\phi: H^{2\ast}(Y, \Zz) \lra
 H^{2\ast}(X, \Zz)$. 
\item
 $Y$ is isomorphic to a fine, two-dimensional moduli space of stable
 sheaves on~$X$.
\end{enumerate}
 \end{Thm}
The non minimal case is covered by the following result of Kawamata.
\begin{Thm}[{\cite[Theorem 1.6]{KawamataKequivalence}}]
Assume that $X$, $Y$ are Fourier-Mukai partners but with $X$ not
minimal.
Then there are only a finite number of possibilities for $Y$ (as in
the minimal case). If $X$ is not isomorphic to a relatively minimal
elliptic rational surface then $X$ and $Y$ are isomorphic. 
\end{Thm}

It remains to classify the auto-equivalences of the derived category
$\cD^b(X)$ for a surface $X$. Orlov solved this problem for an abelian
surface \cite{Orlovabelian} (and more generally for abelian
varieties). Ishii and Uehara \cite{IshiiUehara} solve the problem for the
minimal resolutions of $A_n$-singularities on a surface (so this is a
local result).
The most interesting open case is given by K3-surfaces
although here important progress has recently been made by
Bridgeland~\cite{BridgelandK3,BridgelandStability}. For any $X$
Bridgeland constructs a finite dimensional complex manifold $\Stab(X)$
on which $\Aut(D^b(X))$ acts naturally. Roughly speaking the points of
$\Stab(X)$ correspond to t-structures on $D^b(X)$ together with extra
data defining Harder-Narasimhan filtrations on objects in the heart.
The definition of $\Stab(X)$ was directly inspired by work of Michael
Douglas on stability in string theory \cite{Douglas}.  It seems very
important to obtain a better understanding of the space $\Stab(X)$.

\section{Threefolds and higher dimensional varieties}
\label{ref-6-22}

 If $X$ is a projective smooth
threefolds then just as in the surface case one would like to find 
a unique smooth minimal $X_0$ birationally equivalent to $X$. 
Unfortunately it is well known that this is not possible so
some modifications have to be made. In particular one has
to allow $X_0$ to have some mild singularities, and furthermore
$X_0$ will in general be far from unique.

\emph{Throughout all our varieties are projective}. We say that $X$ is 
\emph{minimal} if $X$ is $\QQ$-Gorenstein and $K_X$ is numerically
effective.
I.e.\ for any curve $C\subset X$ we have $K_X\cdot C\ge 0$. 

A natural category to work in are  varieties with
\emph{terminal} singularities. Recall that a projective variety $X$
has terminal singularities if it is $\QQ$-Gorenstein and for a (any)
resolution $f:Z\rightarrow X$ the discrepancy ($\QQ$-)divisor $K_Z-f^\ast
K_X$ contains every exceptional divisor with strictly
positive coefficients. If $\dim X\le 2$ and $X$ has terminal singularities
then $X$ is smooth. So terminal singularities are indeed very mild. 

If $X$ is a threefold with terminal singularities then
there
exists a map $f:Z\rightarrow X$ which is an isomorphism in codimension one such
that $Z$ terminal, and $\QQ$-factorial \cite[Theorem
6.25]{KollarMori}. Minimal threefolds with 
$\QQ$-factorial terminal singularities are the ``end products'' of 
the three dimensional minimal program. Such minimal models are
however not unique. One has the following classical result by Kollar 
\cite{Kollar}. 
\begin{Thm} Any birational map between minimal
threefolds with $\QQ$-factorial terminal singularities can be
decomposed as a sequence of \emph{flops}.
\end{Thm}
Recall that a flop is a birational map which factors as  $(f^+)^{-1}f$
\[
\xymatrix{%
X\ar[dr]|-f\ar@{.>}[rr] && X^{+}\ar[dl]|-{f+}\\
&W
}
\]
where $f$, $f^+$ are isomorphisms in codimension one
such that $K_X$ and $K_{X^+}$ are $\QQ$-trivial on the fibers of $f$
and $f^+$ respectively and
such that there is a $\QQ$-Cartier divisor $D$ on $X$ with the property
that
$D$ is relatively ample for $f$ and $-D$ is relatively ample for $f^+$.

\begin{Example}\label{ref-6.2-23}
  The easiest (local) example of a flop is the Atiyah flop
  \cite{Reid}: Let $W = \mSpec( \Cc[x,y,z,u]/(xu -yz)$ be the affine
  cone over $\Pp^1 \times \Pp^1$ associated to the line bundle
  $\Oscr_{\PP^1\times \PP^1}(1,1)$. $W$ has an isolated singularity in
  the origin which may be resolved in two different ways $X \lra W
  \longleftarrow X^{+}$ by blowing up the ideals $(x,y)$ and
  $(x,z)$. The varieties $X$ and $X^{+}$ are related by a flop.
\end{Example}

How does one construct a minimal model? Assume that $X$ has
$\QQ$-factorial terminal singularities such that $K_X$ is not
numerically effective. The celebrated cone theorem
\cite{ClemensKollarMori,KollarMori} allows one to construct a map
$f:X\rightarrow W$ with relatively ample $-K_X$ such that one of the
following properties holds \cite[Thm (5.9)]{ClemensKollarMori}
\begin{enumerate}
\item
$\dim X>\dim W$ and $f$ is a $\QQ$-Fano fibration.
\item
$f$ is birational and contracts a divisor.
\item
$f$ is birational and contracts a subvariety of codimension $\ge 2$.
\end{enumerate}
Case 1.\ is what one would get by applying the cone theorem to
$\PP^2$. The result would be the contraction $\PP^2\rightarrow
\mathrm{pt}$. In the case of surfaces 2.\ corresponds to blowing down
exceptional curves. In general the result is again a variety with
terminal singularities and smaller Neron-Severi group.  Case 3.\ 
represents an new phenomenon which only occurs in dimension three and
higher. In this case $W$ may be not be $\QQ$-Gorenstein so one is out
of the category one wants to work in. In order to continue at
this point one introduces a new operation called a \emph{flip}.  A
flip is  a  birational map which factors as $(f^+)^{-1}f$
\[
\xymatrix{%
X\ar[dr]|-f\ar@{.>}[rr] && X^{+}\ar[dl]|-{f+}\\
&W
}
\]
where $f$, $f^+$ are isomorphisms in
codimension one
such that $-K_X$ is relatively ample for $f$, $K_X$ is relatively
ample for $f^+$ and $X^+$ again has $\QQ$-factorial terminal
singularities. The existence of three dimensional flips was settled by
Mori in \cite{Mori}. In higher dimension it is still open.
\begin{Example}\label{ref-6.3-24}
  Let us give an easy example of a (higher dimensional) flip
  generalizing Example \ref{ref-6.2-23}.  Let $W$ be the affine cone over
  $\Pp^m \times \Pp^n$ ($m\le n$) associated to the line bundle
  $\Oscr_{\PP^m\times \PP^n}(1,1)$. $W$ has two canonical resolutions,
  the first one $X$ being given as the total space of the vector
  bundle $\Oscr(1)^{\oplus n}$ over $\PP^m$ and the second one $X^+$
  as the total space of the vector bundle $\Oscr(1)^{\oplus m}$ over
  $\PP^n$. The birational map $X\dashrightarrow X^{+}$ is a flip.
\end{Example}
Following (and slightly generalizing) \cite{BondalOrlovderived} (see
also \cite{KawamataKequivalence}) let us say that a birational map
$X\dashrightarrow X^+$ between $\QQ$-Gorenstein varieties is a
\emph{generalized flip} if there is a commutative diagram with
$\tilde{X}$ smooth
\[
\xymatrix{
&\tilde{X}\ar[ld]|-*+{\pi}\ar[rd]|-*+<0.1ex>{\pi^+} &\\
X\ar@{.>}[rr]&& X^{+}
}
\]
such that $D=\pi^\ast K_X-{\pi^+}^\ast(K_{X^{+}})$ is effective. If
$D=0$ then $X\dashrightarrow X^{+}$ is a generalized flop.

\medskip

Bondal and Orlov \cite{BondalOrlovderived} 
state the following conjecture (see also \cite{KawamataKequivalence}).
\begin{Conjecture} \label{ref-6.4-25}
  For any generalized flip $X\dashrightarrow X^+$ between smooth
  projective varieties there is a full faithful functor $D^b(X^+)\rightarrow
  D^b(X)$. This functor is an equivalence for generalized flops.
\end{Conjecture}
One could think of this conjecture as the foundation for a ``derived
minimal model'' program.

As evidence of the fact that smooth projective varieties related by
a generalized flop are expected to have many properties in common
we recall the following very general result by Batyrev and Kontsevich.
\begin{Thm} 
\label{ref-6.5-26} If $X$ and $X^+$ smooth varieties related by a generalized 
flop then they have the same Hodge numbers.
\end{Thm}
\begin{proof} (see \cite{Batyrev,Craw,Looijenga}) If $X$ and $X^+$ are related by a
  generalized flop then they have the same ``stringy E-function''.
  Since $X$ and $X^+$ are smooth the stringy E-function is equal to
  usual E-function which encodes the Hodge numbers.
\end{proof}
\begin{Remark} The relation 
  between Conjecture \ref{ref-6.4-25} and Theorem \ref{ref-6.5-26} seems rather
  subtle. Indeed a non-trivial Fourier-Mukai transform does not
  usually preserve cohomological degree and hence certainly does not
  preserve the Hodge decomposition.
\end{Remark}

For non-smooth varieties $D^b(X)$ is probably not the correct object
to consider. If $X$ is $\QQ$-Gorenstein then every point $x\in X$ has
some neighborhood $U_x$ such that on $U_x$ there is some positive
number $m_x$ with the property $m_xK_x=0$. Then $K_x$ generates a
cover $\tilde{U}_x$ of $U_x$ on which $\ZZ/m\ZZ$ is acting naturally.
Gluing the local quotient stacks $\tilde{U}_x/(\ZZ/m\ZZ)$ defines a
Deligne-Mumford stack \cite{LaumonMoretBailly} $\Xscr$ birationally
equivalent to $X$.  As usual we write $D^b(\Xscr)$ for
$D^b(\coh(\Xscr))$.  
The following result summarizes what is currently
known in dimension three concering the categories $D^b(\Xscr)$.
\begin{Thm}  
Let $\alpha:X\dashrightarrow X^+$ be a generalized flop
  between threefolds with $\QQ$-factorial terminal
  singularities.
\begin{enumerate}
\item $\alpha$ is a composition of flops.
\item There is a corresponding equivalence $D^b(\Xscr)\rightarrow
  D^b(\Xscr^+)$.
\end{enumerate}
\end{Thm}
In this generality this result was proved by Kawamata in
\cite{KawamataKequivalence}.
The corresponding result in the smooth case
was first proved by Bridgeland in \cite{Bridgelandflops}. By 1)
it is sufficient to consider the case of flops. While trying to
understand Bridgeland's proof the second author produced a mildly
different proof of the result \cite{VandenBerghflops}. Some of the
ingredients in this new proof  
were adapted to the case of stacks by Kawamata. We should also mention
\cite{chen} which uses a different method to extend Bridgeland's result
to singular spaces.

\medskip

Let us give some more comments on flips and flops. 
Flips and flops occur
very naturally in invariant theory \cite{Thaddeus} and 
toric geometry and, as a particular case, for moduli spaces of thin
sincere representations of quivers.

\medskip

Batyrev's
construction of Calabi-Yau varieties \cite{Batyrevmirror} uses toric
geometry, in 
particular toric Fano varieties. Those varieties correspond to
reflexive polytopes. Reflexive polytopes can also be constructed
directly from quivers, however, this class of reflexive polytopes is
very small. For moduli spaces of thin sincere quiver representations
of dimension three all flips are actually flops.  

\begin{Remark} The results above should have consequences for derived
categories of modules over finite dimensional algebras. However, no
example is known of a derived equivalence between a bounded derived
category $\cD^b(A)$ of modules over finite dimensional algebra $A$ and $\cD^b(X)$,
where $X$ admits a flop. The ``closest'' examples to such an equivalence
are the fully faithful functors constructed in \cite{AltmannHille}. If
one allows flips (instead of flops) such equivalences exist, one may
find toric varieties $Y$ with a full strong exceptional sequence of
line bundles. However, for its counterpart $W$ under the flip such
sequences are not known. Strongly related to this problem is a
conjecture of A.~King, that each smooth toric variety admits a full
strong exceptional sequence of line bundles, however, even the
existence of a full exceptional sequence of line bundles is an open
problem (see \cite{Kingtoric} and \cite{AurouxKatzarkovOrlov}). 
\end{Remark}
\section{Non-commutative rings in algebraic geometry}
\label{ref-7-27}
In the previous section we considered mainly Fourier-Mukai transforms
between algebraic varieties. There are also species
of Fourier-Mukai transforms where one of the partners is
non-commutative. In this section we discuss some examples. In contrast
to the previous sections our algebraic varieties will not always be
projective. 

Let $f:X\rightarrow W$ be a projective birational map between Gorenstein
varieties. $f$ is said to be a crepant resolution if $X$ is
smooth and if $f^\ast\omega_W=\omega_X$.  A variant of Conjecture
\ref{ref-6.4-25} is the following:
\begin{Conjecture} \label{ref-7.1-28} Assume that $W$ has Gorenstein singularities and
  that we have two crepant resolutions.
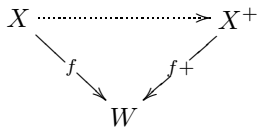
\[
\xymatrix{%
X\ar[dr]|-f\ar@{.>}[rr] && X^{+}\ar[dl]|-{f+}\\
&W
}
\]
Then $X$ and $X^+$ are derived equivalent.
\end{Conjecture}
This conjecture is known in a number of special cases. See the
previous section and \cite{BondalOrlovsemiorthogonal,Bridgelandflops, Namikawa1, KawamataKequivalence}.
There was some initial hope that the derived equivalence between $X$
and $X^+$ would always be induced by $\Oscr_{X\times_{W} X^+}$ but
this turned out to be false for certain so-called ``stratified Mukai-flops''.
See \cite{Namikawa2}.

We will now consider a mild non-commutative situation to which a
similar conjecture applies.  Let $G\subset \Sl_n(\CC)$ be a
finite group and put $W=\CC^n/G$. Write $D^b_G(\CC^n)$ for the
category of $G$ equivariant coherent sheaves on $\CC^n$ and let
$X\rightarrow W$ be a crepant resolution $W$.
\begin{Conjecture} \label{ref-7.2-29} $D^b(X)$ and $D^b_G(\CC^n)$ are derived equivalent.
\end{Conjecture}
If $A$ is the skew group ring $\Oscr(\CC^n)\ast G$ then one may view
$A$ as a \emph{non-commutative crepant resolution} of $\CC^n/G$.  Conjecture
\ref{ref-7.2-29} may be reinterpreted as saying that all commutative crepant
resolutions are derived equivalent to a non-commutative one. So in
that sense it is an obvious generalization of Conjecture \ref{ref-7.1-28}.  A
proper definition of a non-commutative crepant resolution together
with a suitably generalized version of Conjecture \ref{ref-7.2-29} was given
in \cite{VandenBerghCrepant}. An example where this generalized
conjecture applies is \cite{GordonSmith}. A similar but slightly
different conjecture is \cite[Conjecture
5.1]{BondalOrlovderived}.

\medskip

Conjecture \ref{ref-7.2-29} has  now been proved in two cases. First let $X$ be the
irreducible component of the $G$-Hilbert scheme of $\CC^n$ containing
the regular representation. Then we have the celebrated BKR-theorem
\cite{BKR}.
\begin{Thm} Assume that $\dim X\le n+1$ (this holds in particular if
  $n\le 3$). Then $X$ is a crepant resolution of $W$ and $D^b(X)$ is 
equivalent to $D^b_G(\CC^n)$.
\end{Thm}
Note that this theorem, besides establishing the expected derived
equivalence, also produces a specific crepant resolution of $W$. 
For $n=3$ this was done earlier by a case by case analysis
(see \cite{Roan} and the references therein). 

\medskip

Very recently the following result was proved.
\begin{Thm}{\cite{BezrukavnikovKaledin}} Assume that $G$ acts
symplectically on $\CC^n$ (for some arbitrary linear symplectic
form). Then Conjecture \ref{ref-7.2-29} is true.
\end{Thm}
Somewhat surprisingly this result is proved by reduction to
characteristic $p$.
 
\medskip

Let us now discuss a similar but related problem. For a given
scheme $X$ one may want to find algebras $A$ derived equivalent
to $X$. One has the following very general result.
\begin{Thm}[\cite{BondalVandenBergh}, see also \cite{rouquier1}] 
\label{ref-7.5-30} Assume that $X$ is separated. Then
there exists a perfect complex $E$ such that $D(\Qcoh(X))$ is
equivalent to $D(A)$ where $A$ is the DG-algebra $\RHom_{\Oscr_X}(E,E)$.
\end{Thm}
Recall that a perfect complex is one which is locally quasi-isomorphic
to a finite complex of finite rank vector bundles. 

\medskip

In order to replace DG-algebras by real algebras let us say that
a perfect complex $E\in D(\Qcoh(X))$ is \emph{classical tilting} if
it generates $D(\Qcoh(X))$ (in the sense that $\RHom_{\Oscr_X}(E,U)=0$ implies
$U=0$) and $\Hom^i_{\Oscr_X}(E,E)=0$ for $i\neq 0$. One has the
following result.
\begin{Thm} \label{ref-7.6-31} Assume that $X$ is projective over a noetherian affine scheme
  of finite type and assume $E\in D(\Qcoh(X))$ is a classical tilting
  object. Put $A=
\End_{\Oscr_X}(E)$. Then
\begin{enumerate}
\item $\RHom_{\Oscr_X}(E,-)$ induces an equivalence between
  $D(\Qcoh(X))$ and $D(A)$. 
\item This equivalence restricts to an equivalence between
  $D^b(\coh(X))$ and $D^b(\mod(A))$.
\item If $X$ is smooth then $A$ has finite global dimension.
\end{enumerate}
\end{Thm}
\begin{proof}
  1)\ is just a variant on Theorem \ref{ref-7.5-30}.  The inverse functor is
  $-\Lotimes_A E$.  To prove 2) note that the perfect complexes are
  precisely the compact objects (see \cite[Theorem
  3.1.1]{BondalVandenBergh} for a very general version of this
  statement). Hence perfect complexes are preserved under 
$-\Lotimes_A E$. An object $U$ has bounded
  cohomology if and only for any perfect complex $C$ one has
  $\Hom(C,U[n])=0$ for $|n|\gg 0$. Hence objects with bounded
  cohomology are preserved as well. Now let $Z$ be an object in
  $D^b(\mod(A))$. Then it easy to see that $\tau_{\ge n}(Z\Lotimes_A
  E)$ is in $D^b(\coh(X))$ for any $n$.  Since $Z\Lotimes_A E$ has
  bounded cohomology we are done.  To prove 3) note that for any
  $U,V\in\mod(A)$ we have $\Ext^i_A(U,V)$ for $i\gg 0$. Since $A$ has
  finite type this implies that $A$ has finite global dimension.
\end{proof}
Classical tilting objects (and somewhat more generally: ``exceptional
collections'') exist for many classical types of varieties
\cite{BondalOrlovsemiorthogonal}. The following somewhat abstract
result was proved in \cite{VandenBerghflops}.
\begin{Thm} Assume that $f:Y\rightarrow X$ is a projective map between 
varieties, with $X$ affine such that $Rf_\ast\Oscr_Y=\Oscr_X$ and such that
$\dim f^{-1}(x)\le 1$ for all $x\in X$. Then $Y$ has a classical 
tilting object.
\end{Thm} 
This result was inspired by Bridgeland's methods in \cite{Bridgelandflops}.
It applies in particular to resolutions of three-dimensional
Gorenstein terminal singularities. It also has a globalization if
$X$ is quasi-projective instead of affine.


\begin{thebibliography}{10}

\bibitem{AltmannHille}
K.~Altmann and L.~Hille, \emph{Strong exceptional sequences provided by
  quivers}, Algebr. Represent. Theory \textbf{2} (1999), no.~1, 1--17.

\bibitem{AurouxKatzarkovOrlov}
D.~Auroux, L.~Katzarkov, and D.~Orlov, \emph{Mirror symmetry for weighted
  projective planes and their noncommutative deformations}, math.AG/0404281.

\bibitem{Balmer}
P.~Balmer, \emph{Presheaves of triangulated categories and reconstruction of
  schemes}, Math. Ann. \textbf{324} (2002), no.~3, 557--580.

\bibitem{BarthPetersVandeVen}
W.~Barth, C.~Peters, and A.~Van~de Ven, \emph{Compact complex surfaces},
  Ergebnisse der Mathematik und ihrer Grenzgebiete (3) [Results in Mathematics
  and Related Areas (3)], vol.~4, Springer-Verlag, Berlin, 1984.

\bibitem{Batyrevmirror}
V.~V. Batyrev, \emph{Dual polyhedra and mirror symmetry for {C}alabi-{Y}au
  hypersurfaces in toric varieties}, J. Algebraic Geom. \textbf{3} (1994),
  no.~3, 493--535.

\bibitem{Batyrev}
\bysame, \emph{Stringy {H}odge numbers of varieties with {G}orenstein canonical
  singularities}, Integrable systems and algebraic geometry (Kobe/Kyoto, 1997),
  World Sci. Publishing, River Edge, NJ, 1998, pp.~1--32.

\bibitem{Beilinson}
A.~A. Beilinson, \emph{The derived category of coherent sheaves on {${\bf P}\sp
  n$}}, Selecta Math. Soviet. \textbf{3} (1983/84), no.~3, 233--237, Selected
  translations.

\bibitem{BezrukavnikovKaledin}
R.~Bezrukavnikov and D.~Kaledin, \emph{Mc{K}ay equivalence for symplectic
  resolutions of singularities}, math.AG/0401002.

\bibitem{BondalOrlovsemiorthogonal}
A.~Bondal and D.~O. Orlov, \emph{Semi-orthogonal decompositions for algebraic
  varieties}, math.AG 9506012.

\bibitem{BondalOrlovreconstruction}
\bysame, \emph{Reconstruction of a variety from the derived category and groups
  of autoequivalences}, Compositio Math. \textbf{125} (2001), no.~3, 327--344.

\bibitem{BondalOrlovderived}
\bysame, \emph{Derived categories of coherent sheaves}, Proceedings of the
  International Congress of Mathematicians, Vol. II (Beijing, 2002) (Beijing),
  Higher Ed. Press, 2002, pp.~47--56.

\bibitem{BondalVandenBergh}
A.~Bondal and M.~Van~den Bergh, \emph{Generators and representability of
  functors in commutative and noncommutative geometry}, Mosc. Math. J.
  \textbf{3} (2003), no.~1, 1--36, 258.

\bibitem{BondalKapranov}
A.~I. Bondal and M.~M. Kapranov, \emph{Representable functors, {S}erre
  functors, and reconstructions}, Izv. Akad. Nauk SSSR Ser. Mat. \textbf{53}
  (1989), no.~6, 1183--1205, 1337.

\bibitem{borel}
A.~Borel, P.-P. Grivel, B.~Kaup, A.~Haefliger, B.~Malgrange, and F.~Ehlers,
  \emph{Algebraic {$D$}-modules}, Perspectives in Mathematics, vol.~2, Academic
  Press Inc., Boston, MA, 1987.

\bibitem{BrennerButler}
S.~Brenner and M.~C.~R. Butler, \emph{Generalizations of the
  {B}ernstein-{G}el\cprime fand-{P}onomarev reflection functors},
  Representation theory, II (Proc. Second Internat. Conf., Carleton Univ.,
  Ottawa, Ont., 1979), Lecture Notes in Math., vol. 832, Springer, Berlin,
  1980, pp.~103--169.

\bibitem{BridgelandK3}
T.~Bridgeland, \emph{Stability conditions on {K3} surfaces}, math.AG/0307164.

\bibitem{BridgelandStability}
\bysame, \emph{Stability conditions on triangulated categories},
  math.AG/0212237.

\bibitem{Bridgelandequivtriang}
\bysame, \emph{Equivalences of triangulated categories and {F}ourier-{M}ukai
  transforms}, Bull. London Math. Soc. \textbf{31} (1999), no.~1, 25--34.

\bibitem{Bridgelandflops}
\bysame, \emph{Flops and derived categories}, Invent. Math. \textbf{147}
  (2002), no.~3, 613--632.

\bibitem{BKR}
T.~Bridgeland, A.~King, and M.~Reid, \emph{The {M}c{K}ay correspondence as an
  equivalence of derived categories}, J. Amer. Math. Soc. \textbf{14} (2001),
  no.~3, 535--554 (electronic).

\bibitem{BridgeMaciociasurfaces}
T.~Bridgeland and A.~Maciocia, \emph{Complex surfaces with equivalent derived
  categories}, Math. Z. \textbf{236} (2001), no.~4, 677--697.

\bibitem{Caldararu}
A.~Caldararu, \emph{The {M}ukai pairing, {II}: the
  {H}ochschild-{K}ostant-{R}osenberg isomorphism}, math.AG/0308080.

\bibitem{chen}
J.-C. Chen, \emph{Flops and equivalences of derived categories for threefolds
  with only terminal {G}orenstein singularities}, J. Differential Geom.
  \textbf{61} (2002), no.~2, 227--261.

\bibitem{ClemensKollarMori}
H.~Clemens, J.~Koll{\'a}r, and S.~Mori, \emph{Higher-dimensional complex
  geometry}, Ast\'erisque (1988), no.~166, 144 pp. (1989).

\bibitem{Craw}
A.~Craw, \emph{An introduction to motivic integration}, math.AG/9911179.

\bibitem{Deligne}
P.~Deligne, \emph{Cohomologie \'etale}, Springer-Verlag, Berlin, 1977,
  S\'eminaire de G\'eom\'etrie Alg\'ebrique du Bois-Marie SGA 4${1\over 2}$,
  Avec la collaboration de J. F. Boutot, A. Grothendieck, L. Illusie et J. L.
  Verdier, Lecture Notes in Mathematics, Vol. 569.

\bibitem{Douglas}
M.~R. Douglas, \emph{Dirichlet branes, homological mirror symmetry, and
  stability}, Proceedings of the International Congress of Mathematicians, Vol.
  III (Beijing, 2002) (Beijing), Higher Ed. Press, 2002, pp.~395--408.

\bibitem{gabriel}
P.~Gabriel, \emph{Des cat\'egories ab\'eliennes}, Bull. Soc. Math. France
  \textbf{90} (1962), 323--448.

\bibitem{GeigleLenzing}
W.~Geigle and H.~Lenzing, \emph{A class of weighted projective curves arising
  in representation theory of finite-dimensional algebras}, Singularities,
  representation of algebras, and vector bundles (Lambrecht, 1985), Lecture
  Notes in Math., vol. 1273, Springer, Berlin, 1987, pp.~265--297.

\bibitem{GordonSmith}
I.~Gordon and S.~P. Smith, \emph{Representations of symplectic reflection
  algebras and resolutions of deformations of symplectic quotient
  singularities}, math.RT/0310187.

\bibitem{GriffithsHarris}
P.~Griffiths and J.~Harris, \emph{Principles of algebraic geometry}, Wiley
  Classics Library, John Wiley \& Sons Inc., New York, 1994, Reprint of the
  1978 original.

\bibitem{Happel}
D.~Happel, \emph{Triangulated categories in the representation theory of
  finite-dimensional algebras}, London Mathematical Society Lecture Note
  Series, vol. 119, Cambridge University Press, Cambridge, 1988.

\bibitem{RingelHappelTubular}
D.~Happel and C.~M. Ringel, \emph{The derived category of a tubular algebra},
  Representation theory, I (Ottawa, Ont., 1984), Lecture Notes in Math., vol.
  1177, Springer, Berlin, 1986, pp.~156--180.

\bibitem{Hilton}
P.~Hilton, \emph{General cohomology theory and {$K$}-theory}, Course given at
  the University of S\~ao Paulo in the summer of 1968 under the auspices of the
  Instituto de Pesquisas Matem\'aticas, Universidade de S\~ao Paulo. London
  Mathematical Society Lecture Note Series, vol.~1, Cambridge University Press,
  London, 1971.

\bibitem{Huybrechts}
D.~Huybrechts, \emph{Fourier-{M}ukai transforms in algebraic geometry},
  \texttt{http://www.institut.math.jussieu.fr/\~{}huybrech/FM.ps}, preliminary
  lecture notes.

\bibitem{IshiiUehara}
A.~Ishii and H.~Uehara, \emph{{Autoequivalences of derived categories on the
  minimal resolutions of $A_n$-singularities on surfaces}}, math.AG/0409151.

\bibitem{Karoubi}
M.~Karoubi, \emph{{$K$}-theory}, Springer-Verlag, Berlin, 1978, An
  introduction, Grundlehren der Mathematischen Wissenschaften, Band 226.

\bibitem{kashiwara}
M.~Kashiwara, \emph{The {R}iemann-{H}ilbert problem for holonomic systems},
  Publ. Res. Inst. Math. Sci. \textbf{20} (1984), no.~2, 319--365.

\bibitem{KawamataKequivalence}
Y.~Kawamata, \emph{{$D$}-equivalence and {$K$}-equivalence}, J. Differential
  Geom. \textbf{61} (2002), no.~1, 147--171.

\bibitem{Kawamata}
\bysame, \emph{Equivalences of derived categories of sheaves on smooth stacks},
  Amer. J. Math. \textbf{126} (2004), no.~5, 1057--1083.

\bibitem{KS}
M.~Khovanov and P.~Seidel, \emph{Quivers, {F}loer cohomology, and braid group
  actions}, J. Amer. Math. Soc. \textbf{15} (2002), no.~1, 203--271
  (electronic).

\bibitem{Kingtoric}
A.~King, \emph{Tilting bundles on some rational surfaces}, preprint
  \texttt{http://www.maths.bath.ac.uk/\~{}masadk/papers/}.

\bibitem{Kollar}
J.~Koll{\'a}r, \emph{Flops}, Nagoya Math. J. \textbf{113} (1989), 15--36.

\bibitem{KollarMori}
J.~Koll{\'a}r and S.~Mori, \emph{Birational geometry of algebraic varieties},
  Cambridge Tracts in Mathematics, vol. 134, Cambridge University Press,
  Cambridge, 1998, With the collaboration of C. H. Clemens and A. Corti,
  Translated from the 1998 Japanese original.

\bibitem{Kontsevichhomolmirror}
M.~Kontsevich, \emph{Homological algebra of mirror symmetry}, Proceedings of
  the International Congress of Mathematicians, Vol.\ 1, 2 (Z\"urich, 1994)
  (Basel), Birkh\"auser, 1995, pp.~120--139.

\bibitem{LaumonMoretBailly}
G.~Laumon and L.~Moret-Bailly, \emph{Champs alg\'ebriques}, Ergebnisse der
  Mathematik und ihrer Grenzgebiete. 3. Folge. A Series of Modern Surveys in
  Mathematics [Results in Mathematics and Related Areas. 3rd Series. A Series
  of Modern Surveys in Mathematics], vol.~39, Springer-Verlag, Berlin, 2000.

\bibitem{MeltzerLenzingwtprojlinegenus1}
H.~Lenzing and H.~Meltzer, \emph{Sheaves on a weighted projective line of genus
  one and representations of a tubular algebra}, Proceedings of the Sixth
  International Conference on Representations of Algebras (Ottawa, ON, 1992)
  (Ottawa, ON), Carleton-Ottawa Math. Lecture Note Ser., vol.~14, Carleton
  Univ., 1992, p.~25.

\bibitem{Looijenga}
E.~Looijenga, \emph{Motivic measures}, Ast\'erisque (2002), no.~276, 267--297,
  S\'eminaire Bourbaki, Vol.\ 1999/2000.

\bibitem{mebkhout2}
Z.~Mebkhout, \emph{Une autre \'equivalence de cat\'egories}, Compositio Math.
  \textbf{51} (1984), no.~1, 63--88.

\bibitem{mebkhout1}
\bysame, \emph{Une \'equivalence de cat\'egories}, Compositio Math. \textbf{51}
  (1984), no.~1, 51--62.

\bibitem{MiyachiYekutieli}
J.-i. Miyachi and A.~Yekutieli, \emph{Derived {P}icard groups of
  finite-dimensional hereditary algebras}, Compositio Math. \textbf{129}
  (2001), no.~3, 341--368.

\bibitem{Mori}
S.~Mori, \emph{Flip theorem and the existence of minimal models for
  {$3$}-folds}, J. Amer. Math. Soc. \textbf{1} (1988), no.~1, 117--253.

\bibitem{Mukai}
S.~Mukai, \emph{Duality between {$D(X)$} and {$D(\hat X)$} with its application
  to {P}icard sheaves}, Nagoya Math. J. \textbf{81} (1981), 153--175.

\bibitem{MukaiK3}
\bysame, \emph{On the moduli space of bundles on {$K3$} surfaces. {I}}, Vector
  bundles on algebraic varieties (Bombay, 1984), Tata Inst. Fund. Res. Stud.
  Math., vol.~11, Tata Inst. Fund. Res., Bombay, 1987, pp.~341--413.

\bibitem{Namikawa1}
Y.~Namikawa, \emph{Mukai flops and derived categories}, J. Reine Angew. Math.
  \textbf{560} (2003), 65--76.

\bibitem{Namikawa2}
\bysame, \emph{Mukai flops and derived categories. {II}}, Algebraic structures
  and moduli spaces, CRM Proc. Lecture Notes, vol.~38, Amer. Math. Soc.,
  Providence, RI, 2004, pp.~149--175.

\bibitem{OrlovFMtransf}
D.~O. Orlov, \emph{Equivalences of derived categories and {$K3$} surfaces}, J.
  Math. Sci. (New York) \textbf{84} (1997), no.~5, 1361--1381, Algebraic
  geometry, 7.

\bibitem{Orlovabelian}
\bysame, \emph{Derived categories of coherent sheaves on abelian varieties and
  equivalences between them}, Izv. Ross. Akad. Nauk Ser. Mat. \textbf{66}
  (2002), no.~3, 131--158.

\bibitem{OrlovRussianMathSurveys}
\bysame, \emph{Derived categories of coherent sheaves and equivalences between
  them}, Uspekhi Mat. Nauk \textbf{58} (2003), no.~3, 89--172.

\bibitem{Reid}
M.~Reid, \emph{What is a flip?}, colloquium talk at Utah 1992.

\bibitem{Rickard1}
J.~Rickard, \emph{Morita theory for derived categories}, J. London Math. Soc.
  (2) \textbf{39} (1989), no.~3, 436--456.

\bibitem{Rickard2}
\bysame, \emph{Derived equivalences as derived functors}, J. London Math. Soc.
  (2) \textbf{43} (1991), no.~1, 37--48.

\bibitem{Ringeltame}
C.~M. Ringel, \emph{Tame algebras and integral quadratic forms}, Lecture Notes
  in Mathematics, vol. 1099, Springer-Verlag, Berlin, 1984.

\bibitem{Roan}
S.~Roan, \emph{Minimal resolutions of {G}orenstein orbifolds in dimension
  three}, Topology \textbf{35} (1996), no.~2, 489--508.

\bibitem{rosenberg}
A.~L. Rosenberg, \emph{The spectrum of abelian categories and reconstruction of
  schemes}, Rings, Hopf algebras, and Brauer groups (Antwerp/Brussels, 1996),
  Lecture Notes in Pure and Appl. Math., vol. 197, Dekker, New York, 1998,
  pp.~257--274.

\bibitem{rouquier1}
R.~Rouquier, \emph{Dimensions of triangulated categories}, math.CT/0310134.

\bibitem{rouquier}
\bysame, \emph{Cat\'egories d\'eriv\'ees et g\'eom\'etrie alg\'ebrique},
  \texttt{http://www.math.jussieu.fr/\~{}rouquier/preprints/luminy.pdf}, 2004,
  preprint.

\bibitem{RZ}
R.~Rouquier and A.~Zimmermann, \emph{Picard groups for derived module
  categories}, Proc. London Math. Soc. (3) \textbf{87} (2003), no.~1, 197--225.

\bibitem{Sato}
M.~Sato, T.~Kawai, and M.~Kashiwara, \emph{Microfunctions and
  pseudo-differential equations}, Hyperfunctions and pseudo-differential
  equations (Proc. Conf., Katata, 1971; dedicated to the memory of Andr\'e
  Martineau), Springer, Berlin, 1973, pp.~265--529. Lecture Notes in Math.,
  Vol. 287.

\bibitem{Schiffmann}
O.~Schiffmann, \emph{Noncommutative projective curves and quantum loop
  algebras.}, math.QA 0205267.

\bibitem{Seidel}
P.~Seidel, \emph{Homological mirror symmetry for the quartic surface},
  math.SG/0310414.

\bibitem{Seidel1}
\bysame, \emph{Lectures on four-dimensional {D}ehn twists}, math.SG/0309012.

\bibitem{SeidelThomas}
P.~Seidel and R.~Thomas, \emph{Braid group actions on derived categories of
  coherent sheaves}, Duke Math. J. \textbf{108} (2001), no.~1, 37--108.

\bibitem{Thaddeus}
M.~Thaddeus, \emph{Geometric invariant theory and flips}, J. Amer. Math. Soc.
  \textbf{9} (1996), no.~3, 691--723.

\bibitem{Ueno}
K.~Ueno, \emph{Classification theory of algebraic varieties and compact complex
  spaces}, Springer-Verlag, Berlin, 1975, Notes written in collaboration with
  P. Cherenack, Lecture Notes in Mathematics, Vol. 439.

\bibitem{VandenBerghCrepant}
M.~Van~den Bergh, \emph{Non-commutative crepant resolutions}, The legacy of
  Niels Henrik Abel, Springer, Berlin, 2004, pp.~749--770.

\bibitem{VandenBerghflops}
\bysame, \emph{Three-dimensional flops and noncommutative rings}, Duke Math. J.
  \textbf{122} (2004), no.~3, 423--455.

\bibitem{Verdier}
J.-L. Verdier, \emph{Des cat\'egories d\'eriv\'ees des cat\'egories
  ab\'eliennes}, Ast\'erisque (1996), no.~239, xii+253 pp. (1997), With a
  preface by Luc Illusie, Edited and with a note by Georges Maltsiniotis.

\end{thebibliography}
\def\cprime{$'$}
\providecommand{\bysame}{\leavevmode\hbox to3em{\hrulefill}\thinspace}
\providecommand{\MR}{\relax\ifhmode\unskip\space\fi MR }
\providecommand{\MRhref}[2]{%
  \href{http://www.ams.org/mathscinet-getitem?mr=#1}{#2}
}
\providecommand{\href}[2]{#2}

\medskip

{{\small Lutz Hille}\\
{\small Mathematisches Seminar}\\
{\small Universit\"at Hamburg}\\
{\small D-20 146 Hamburg}\\
{\small Germany}\\
{\small E-mail: hille@math.uni-hamburg.de}\\
{\small http://www.math.uni-hamburg.de/home/hille/}}

\medskip

{{\small Michel Van den Bergh}\\
{\small Departement WNI}\\
{\small Limburgs Universitair Centrum}\\
{\small Universitaire Campus}\\
{\small 3590 Diepenbeek}\\
{\small Belgium\\}
{\small E-mail: vdbergh@luc.ac.be}\\
{\small http://alpha.luc.ac.be/Research/Algebra/Members/michel\_id.html}}

\end{document}